\documentclass[a4paper]{article}
\usepackage{amsfonts}
\usepackage[english]{babel}
\usepackage{amsmath,mathrsfs,amssymb,amsthm,epsfig,color,graphicx}
\usepackage[latin1]{inputenc}

\numberwithin{equation}{section}
\newtheorem{theorem}{Theorem}[section]
\newtheorem{proposition}[theorem]{Proposition}

\newtheorem{remark}[theorem]{Remark}

\definecolor{Red}{cmyk}{0,1,1,0}

\definecolor{Blue}{cmyk}{1,1,0,0}

\def\fin { \vskip 0pt \hfill \hbox{\vrule height 5pt width 5pt depth 0pt} \vskip 12pt}

\begin{document}

\title{On nonlinear Schrödinger equations with random potentials: existence
and probabilistic properties}
\author{{Leandro Cioletti} \\
{\small Departamento de Matemática, UnB, 70910-900 Brasília, Brazil.}\\
{\small \texttt{E-mail:leandro.mat@gmail.com}}\vspace{.5cm}\\
{Lucas C. F. Ferreira}\\
{\small Universidade Estadual de Campinas, IMECC - Departamento de Matemá%
tica,} \\
{\small {Rua Sérgio Buarque de Holanda, 651, CEP 13083-859, Campinas-SP,
Brazil.}}\\
{\small \texttt{E-mail:lcff@ime.unicamp.br}} \vspace{.5cm}\\
{Marcelo Furtado}\\
{\small Departamento de Matemática, UnB, 70910-900 Brasília, Brazil. }\\
{\small \texttt{E-mail:mfurtado@unb.br}}}
\date{}
\maketitle

\begin{abstract}
In this paper we are concerned with nonlinear Schrödinger equations with
random potentials. Our class includes continuum and discrete potentials.
Conditions on the potential $V_{\omega }$ are found for existence of
solutions almost sure $\omega $. We study probabilistic properties like
central limit theorem and law of larger numbers for the obtained solutions
by independent ensembles. We also give estimates on the expected value for
the $L^{\infty }$-norm of the solution showing how it depends on the size of
the potential.
\end{abstract}

\medskip \noindent {\small {\textbf{AMS 2000 subject classification:} 47B80,
60H25, 35J60, 35R60, 82B44, 47H10 }} \smallskip \noindent {\small {\textbf{%
Keywords:} Random potentials; Random nonlinear equations;
Schrödinger operators}}

\section{Introduction}

\ \ \ \ A class of models that appears naturally in a wide number of
phenomena are the random differential equations. This occurs because
randomness is a powerful tool and concept to control complex systems
involving a large number of variables and particles. The basic idea is
describe complex systems by means of their statistical properties. Another
kind of phenomena are those governed by quantum mechanics and uncertainty
principle. In this direction, we have Schrödinger equations, and their
random versions, which are core in the study of condensed matter.

In this paper we are concerned with a random version of the
nonlinear Schrödinger equation
\begin{equation}
ih\dfrac{\partial \psi }{\partial t}=-h^{2}\Delta \psi +V(x)\psi -|\psi
|^{p-1}\psi ,~~x\in \mathbb{R}^{n},  \label{sch1}
\end{equation}%
where $t\in \mathbb{R}$, $n\geq 3$, $p>1$, $h$ is the Planck constant and $i$
is the imaginary unit. When looking for standing wave solutions, namely
those which have the special form $\psi (x,t):=e^{-i\frac{E}{h}t}u(x)$, with
$E\in \mathbb{R}$, we are leading to solve the following stationary equation
\begin{equation*}
-\Delta u+V(x)u=|u|^{p-1}u,~~~x\in \mathbb{R}^{N}.
\end{equation*}
From the physical viewpoint, the function $V$ is the potential
energy, and therefore the force acting on the system is given by
$F(x)=-\nabla V(x)$. In the deterministic case, there are many
papers concerning existence, multiplicity and qualitative
properties for the solution of the above equation (see
\cite{Rab,delFel,AmbBadCin,AmbMalSec} and references therein).

The main interest of this paper is to study situations where the potential $%
V $ is not deterministic. Worth to mention that during the last thirty
years, random Schrödinger operators, which originated in condensed matter
physics, have been studied intensively by physicists and mathematicians. The
theory is at the crossroads of a number of mathematical fields: the theory
of operators, partial differential equations, the theory of probabilities
and also stochastic process. This paper aims to prove the existence and
probabilistic properties of bounded solutions for the random equation
\begin{equation}
\left\{
\begin{array}{rcll}
-\Delta u+V_{\omega }(x)u & = & b(x)u|u|^{p-1}+g(x), & \text{if}\ x\in U; \\%
[0.15cm]
u & = & 0, & \text{if}\ x\in \partial U,%
\end{array}%
\right.  \label{eqdif-est1}
\end{equation}%
where $V_{\omega }$ is a random variable, $U\subset \mathbb{R}^{n}$ is a
bounded domain and the terms $b,\,g\in L^{\infty }(U)$ are deterministic. In
fact, the boundedness of $U$ is not essential and could be circumvented by
working in weighted $L^{\infty }$-spaces or Lebesgue spaces $L^{s}(\mathbb{R}%
^{n})$ with $s\neq \infty $ (see \cite{Fer1,Fer2}). However, here this
condition will simplify matters a bit. The random potential $V_{\omega }$ is
constructed via a convolution with a realization of a random variable valued
in the finite random measure space. Precisely, given a continuous function $%
f:\mathbb{R}^{N}\rightarrow \mathbb{R}$ we consider
\begin{equation}
V_{\omega }(x):=\int_{U}f(x-y)\,d{\mu _{\omega }}(y)  \label{aux-pot-int}
\end{equation}%
where $\mu _{\omega }$ is a $\mathcal{M}(U)$-valued random variable.

We present here some examples of (\ref{aux-pot-int}) that have been treated
in the literature (see e.g. the review \cite{Kir}). We first consider a
model of an unordered alloy, that is, a mixture of several materials with
atoms located at lattice positions. If we assume that the type of atom at
the lattice $i\in \mathbb{Z}^{n}$ is random we are leading to consider the
following type of potential
\begin{equation}
V_{\omega }(x)=\sum_{i\in \mathbb{Z}^{n}}q_{i}(\omega )f(x-i),  \label{mod1}
\end{equation}%
where the random variables $q_{i}$ describe the charge of the atom at the
position $i$ of the lattice. Other example can be obtained if we consider
materials like glass or rubber, where the position of the atoms of the
material are located at random points $\eta _{i}$ in space. By normalizing
the charge of the atoms, the suggested potential is formally
\begin{equation}
V_{\omega }(x)=\sum_{i\in \mathbb{Z}^{n}}f(x-\eta _{i}(\omega )),
\label{mod2}
\end{equation}%
where the $\eta _{i}(\omega )$ are random variables which localize the atoms
in the spaces.

The class of potentials allowed here is sufficient large to consider many
known models. For example, the case of glass considered in \eqref{mod2} can
be obtained if we take the random point measure $\mu _{\omega
}=\sum_{i}\delta _{\eta _{i}(\omega )}$. Actually, for this choice of the
measure we have that
\begin{equation}
\sum_{i\in \mathbb{Z}^{n}\cap U}f(x-\eta _{i}(\omega ))=\int_{U}f(x-\eta
)d\mu _{\omega }(\eta ).  \label{mod3}
\end{equation}%
Also, a combination of potentials like (\ref{mod1}) and (\ref{mod2}), namely
$\Sigma _{_{i\in \mathbb{Z}^{n}\cap U}}q_{i}(\omega )f(x-\eta _{i}(\omega ))$
(see \cite{ComHis}), is also covered by (\ref{aux-pot-int}) with $\mu
_{\omega }=\Sigma _{_{i\in \mathbb{Z}^{n}\cap U}}q_{i}(\omega )\delta _{\eta
_{i}(\omega )}$. It is not difficult to see that we can also consider other
models like, e.g., the Poisson model (see \cite{Kir} for more examples).

The models (\ref{mod1}) and (\ref{mod2}) correspond to discrete measures $%
\mu _{\omega }$ and results for them about localization, spectral properties
or decays can be found in \cite{BourKen,ComHis,GerKleinSchen,Kir,Saf}. For
Schrödinger equations defined in a lattice, that is $x\in \mathbb{Z}^{n}$,
we refer the reader to \cite{Bourgain,BourWan}. Considering a random
time-dependent potential for (\ref{sch1}), the authors of \cite{GuiTomLen}
studied asymptotic behavior of solutions by showing convergence for
stochastic Gaussian limits when the two-point correlation function of the
potential is rapidly decaying. Still for time random potentials, scaling
limits for parabolic waves in random media were investigated in \cite{Fan}.
Despite important progress in the last years, there is still a lack of
result for random equations, including Schrödinger ones (see \cite{Carmona}%
), mainly with respect to the continuum case which seems to be
harder-to-treating. Another type of random equations are the parabolic ones,
for which we refer the works \cite{JosAli,DawKou} and their references.

In this paper we find conditions on the potential $V_{\omega }$ for the
nonlinear equation (\ref{eqdif-est1}) having solutions almost sure $\omega .$
The solution are understood in an integral sense coming from Green
functions. From Theorem \ref{teo-moment} we see how the expected value of
the $L^{\infty }$-norm of solutions depends on the size of potential. Our
results also cover continuum random potential like, among others, the
examples given in Remark \ref{Rem1} and Theorem \ref{borel-canteli}.
Moreover, we study probabilistic properties like central limit theorem and
law of larger numbers for the obtained solutions by independent ensembles.
It is worthwhile to mention that, when dealing with the random variable $%
\omega \mapsto u(x,\omega)$ which maps an element of $\Omega$ in the
solution of (\ref{eqdif-est1}) associated with the random potential $%
V_{\omega}$, we need to extend some known concepts of real random variables
for that taking values in a more general Banach space. We refer to Section 2
for more details.

As a further comment, we observe that the random potentials considered in
this paper are built from a very general probability space. In this setting
does not always make sense to ask what is the probability that the problem %
\eqref{eqdif-est1} has an unique solution in $L^{\infty }(U)$. In order to
give some sense to this question we should restrict ourself to probability
spaces $(\Omega ,\mathcal{F},\mathbb{P})$ and random potentials $V$ where
the set
\begin{equation*}
\{\omega \in \Omega :\ \text{the problem \eqref{eqdif-est1} has a unique
solution in}\ L^{\infty }(U)\}
\end{equation*}%
is an event (measurable). Working in such probability spaces Theorem \ref%
{teo-prob} give us immediately a lower bound for the probability that the
non-linear problem \eqref{eqdif-est1} has a unique solution.

The manuscript is organized as follows. In the next section, we introduce
some notations, basic definitions and give some properties for an integral
operator associated with the random potential $V_{\omega }.$ The results are
stated and proved in Section 3.

\section{Preliminaries and notation}

\ \ \ \ Throughout this paper $(\Omega ,\mathcal{F},\mathbb{P})$ denotes a
given complete probability space. If $(E,\mathcal{E})$ is a measurable
space, any $(\mathcal{F},\mathcal{E})$-measurable function $X:\Omega
\rightarrow E$ will be called a $E$-valued random variable. We use the
abbreviation \textit{a.s.} for \textit{almost surely} or \textit{almost sure}.

Let $U\subset \mathbb{R}^{n}$ be a bounded domain. We adopt the standard
notation $\mathcal{M}(U)$ to denote the set of all Random measures over $U$
having finite variation and we call $\mathscr{B}(\mathcal{M}(U))$ the $%
\sigma $-algebra of the borelians of $\mathcal{M}(U)$ generated by the total
variation norm. The space of all bounded continuous real-valued functions
defined on $U$ will be denoted by $BC(U)$. Since $BC(U)$ is a metric space
with the supremum norm, when we refer to a $BC(U)$-valued random variable,
the $\sigma $-algebra we are considering is always the one generated by the
borelians. Similarly to a $\mathcal{X}$-valued Borel random variable $%
X:\Omega \rightarrow \mathcal{X},$ where $\mathcal{X}$ is an arbitrary
metric space.

The random potentials considered here are the $BC(U)$-valued random
variables defined as follows. Take any random variable $X:\Omega \rightarrow
\mathcal{M}(U)$ (which is simply a random measure in $\mathcal{M}(U)$) and a
fixed function $f\in BC(\mathbb{R}^{n})$. Then, for $\mu _{\omega }=X(\omega
)$, the function $V:\Omega \rightarrow BC(U)$ defined by
\begin{equation*}
V_{\omega }(x):=\int_{U}f(x-y)\,d\mu _{\omega }(y),~~~x\in U,
\end{equation*}%
is a $BC(U)$-valued random variable that will be called a random potential.
To see that $V$ is a well-defined $BC(U)$-valued random variable, is enough
to consider the mapping $T_{f}:\mathcal{M}(U)\rightarrow BC(U)$ given by
\begin{equation*}
T_{f}(\mu )(x)=\int_{U}f(x-y)\,d\mu (y),~~~x\in U,
\end{equation*}%
and to observe that $V=T_{f}\circ X$. In fact, if we denote by $%
\rule[-0.8mm]{0.4mm}{3.4mm}\,\mu \,\rule[-0.8mm]{0.4mm}{3.4mm}$ the total
variation of the measure $\mu $, the inequality
\begin{equation}
\Vert T_{f}(\mu )\Vert _{\infty }:=\sup_{x\in U}|T_{f}(\mu )(x)|\leq \left(
\sup_{x\in \mathbb{R}^{n}}\left\vert f(x)\right\vert \right) \,%
\rule[-0.8mm]{0.4mm}{3.4mm}\,\mu \,\rule[-0.8mm]{0.4mm}{3.4mm}
\label{est-medida}
\end{equation}%
implies that $T_{f}$ is a continuous and Borel measurable function. Since $V$
is a composition of two Borel measurable functions, $V$ is a $BC(U)$-valued
random variable.

As usual, if $(U,\mathscr{B},\mu )$ is a measure space, we define
\begin{equation*}
\Vert f\Vert _{L^{\infty }(U,d\mu )}=\inf \left\{ a\geq 0:\mu
(\{x:|f(x)|>a\})=0\right\}
\end{equation*}%
and the space $L^{\infty }(U,\mathscr{B}(U),\mu )$ as being the set
\begin{equation*}
\{f:U\rightarrow \mathbb{R}:f\ \text{is Borel measurable and}\ \Vert f\Vert
_{L^{\infty }(U,d\mu )}<\infty \}.
\end{equation*}%
When $d\mu =dx$ is the Lebesgue measure in $U\subset \mathbb{R}^{n},$ we
simply denote $L^{\infty }(U)=L^{\infty }(U,\mathscr{B}(U),dx)$. Although we
are assuming that $f\in BC(\mathbb{R}^{n})$, most of the results presented
here are also valid if we suppose only the weaker condition $f\in \cap _{\mu
\in \mathcal{M}(U-U)}L^{\infty }(U-U,\mathscr{B}(U-U),\mu )$.

In order to state some convergence results obtained in this paper we need to
use the notion of Bochner integrals. Let $(\mathcal{X},\Vert \cdot \Vert _{%
\mathcal{X}})$ be a Banach space and $(\Omega ,\mathcal{F},\mathbb{P})$ be a
probability space. If $X:\Omega \rightarrow \mathcal{X}$ is a $\mathcal{X}$%
-valued Borel random variable such that $X=Y$ a.s. in $\Omega ,$ where $%
Y:\Omega \rightarrow \mathcal{X}$ is a $\mathcal{X}$-valued Borel random
variable with $Y(\Omega )\subset \mathcal{X}$ separable, and
\begin{equation*}
\int_{\Omega }\Vert X(\omega )\Vert _{\mathcal{X}}\,d\mathbb{P}(\omega
)<\infty ,
\end{equation*}%
then there exist a unique element $\mathbb{E}[X]\in \mathcal{X}$ with the
property
\begin{equation*}
\ell (\mathbb{E}[X])=\int_{\Omega }\ell (X(\omega ))\,d\mathbb{P}(\omega )
\end{equation*}%
for all $\ell \in \mathcal{X}^{\ast }$, where $\mathcal{X}^{\ast }$ is the
dual of $\mathcal{X}$. Following the standard notation we write
\begin{equation*}
\mathbb{E}[X]=\int_{\Omega }X(\omega )\,d\mathbb{P}(\omega ).
\end{equation*}%
We call $\mathbb{E}[X]$ the Bochner integral of $X$ with respect to $\mathbb{%
P}$. More details about the existence and some properties of this integral
can be found in \cite{HilPhi,Par}. For these $\mathcal{X}$-valued random
variables we define the convergence in probability similarly to the
real-valued case, that is, if $\{X_{j}\}$ is a sequence of $\mathcal{X}$%
-valued random variable we say that $X_{j}$ converges to a $\mathcal{X}$%
-valued random variable $X$ in probability if for all $\varepsilon >0$, we
have
\begin{equation}
\lim_{j\rightarrow \infty }\mathbb{P}(\{\omega \in \Omega :\Vert
X_{j}(\omega )-X(\omega )\Vert _{\mathcal{X}}\geq \varepsilon \})=0.
\label{def-conv-prob}
\end{equation}%
When $X$ is real-valued random variable, we use the usual notation and
denote the expected value of $X$ and its variance by
\begin{equation*}
\mathbb{E}[X]:=\int_{\Omega }X(\omega )\,d\mathbb{P}(\omega )\text{ \ \ and
\ Var}\,X:=\mathbb{E}[(\mathbb{E}[X]-X)^{2}],
\end{equation*}%
respectively. For the both senses of expectation presented above we also use
the notation
\begin{equation*}
\mathbb{E}_{A}[X]=\int_{A}X(\omega )\,d\mathbb{P}(\omega ),
\end{equation*}%
whenever $A\subset \Omega $ is measurable and the right-hand-side of the
expression makes sense.

Let $X$ and $Y$ be two $E$-valued random variable in the same probability
space. We say that they are identically distributed if for all $A\in
\mathcal{E}$ we have $\mathbb{P}(X^{-1}(A))=\mathbb{P}(Y^{-1}(A))$. Now we
introduce the notion of independence. Given a finite set of random variables
$X_{1},\ldots X_{j}$ we say they are independent if for all $A_{i}\in
\mathcal{E},1\leq i\leq j$, we have
\begin{equation*}
\mathbb{P}(\cap _{i=1}^{j}X_{i}\in A_{i})=\prod_{i=1}^{j}\mathbb{P}(X_{i}\in
A_{i}).
\end{equation*}%
Finally a sequence of random variables $\{X_{1},X_{2}\ldots \}$ is said
independent if all finite collection of this sequence form a set of
independent random variables. If $X_{1},X_{2},\ldots $ is a sequence of
independent and identically distributed random variables we say that $%
X_{1},X_{2},\ldots $ are i.i.d. random variables.

\section{Main results and proofs}

\ \ \ \ Let $G$ be the Green function of the laplacian operator $-\Delta $
in the bounded domain $U\subset \mathbb{R}^{n}$ with $n\geq 3.$ It is known
that, for all $x,\,y \in U$, there holds
\begin{equation*}
0\leq G(x,y)\leq \frac{1}{n\alpha _{n}(n-2)}\frac{1}{|x-y|^{n-2}},~
\label{est-green}
\end{equation*}
where $\alpha _{n}$ stands for the volume of the unit ball in $\mathbb{R}%
^{n} $. Hence, if we denote by $d_{U}$ the the diameter of $U$, namely
\begin{equation*}
d_{U}:=\sup_{x_{1},\,x_{2}\in U}{|}x_{1}-x_{2}{|},
\end{equation*}%
and $B_{d_{U}}(x)=\{x\in \mathbb{R}^{n};\left\vert x\right\vert <d_{U}\},$ a
straightforward calculation provides

\begin{equation}  \label{est-int-green}
\begin{array}{lcl}
\displaystyle\int_{U}G(x,y)dy & \leq & \dfrac{1}{n\alpha _{n}(n-2)}%
\displaystyle\int_{B_{d_{U}}(x)}\dfrac{1}{|x-y|^{n-2}}dy \vspace{0.2cm} \\
& = & \dfrac{1}{n\alpha _{n}(n-2)}\dfrac{n\alpha _{n}d_{U}^{2}}{2}=\dfrac{%
d_{U}^{2}}{2(n-2)},%
\end{array}%
\end{equation}
for all $x \in U$. From now on we write only $l_{0}=l_{0}(n,U)$ to denote
the following quantity
\begin{equation}
l_{0}:=\dfrac{d_{U}^{2}}{2(n-2)}.  \label{def-l0}
\end{equation}

Inequality \eqref{est-int-green} implies that is well defined the map $%
H:L^{\infty }(U)\rightarrow L^{\infty }(U)$ given by
\begin{equation*}
H(\varphi )(x):=\int_{U}G(x,y)\varphi (y)dy,~~~x \in U.
\end{equation*}%
More specifically, for any $\varphi \in L^{\infty }(U )$, there holds
\begin{equation*}
|H(\varphi )(x)|\leq \int_{U}G(x,y)|\varphi (y)|dy\leq \Vert \varphi \Vert
_{\infty }\int_{U}G(x,y)dy
\end{equation*}%
and therefore
\begin{equation}
\Vert H(\varphi )\Vert _{\infty }\leq l_{0}\Vert \varphi \Vert _{\infty }.
\label{est-h}
\end{equation}

Standard calculations show that the problem \eqref{eqdif-est1} is formally
equivalent to the integral equation
\begin{equation}
u(x)=H(g)+H(V_{\omega }u)+H(bu|u|^{p-1}).  \label{eq-integral}
\end{equation}%
In what follows we make suitable estimates on the terms of the integral
equation in order to be able to apply a fixed point argument. We first set $%
\mathcal{X}:=L^{\infty }(U)$ and define, for any fixed $\omega \in \Omega $,
the linear function $T:\mathcal{X}\rightarrow \mathcal{X}$ by
\begin{equation*}
T(u):=H(V_{\omega }u),~~~\forall \,u\in \mathcal{X}.
\end{equation*}%
It follows from \eqref{est-h} and \eqref{est-medida} that, for any $u\in
\mathcal{X}$, there holds
\begin{equation}
\Vert T(u)\Vert _{\infty }\leq l_{0}\Vert V_{\omega }u\Vert _{\infty }\leq
l_{0}\Vert f\Vert _{\infty }\rule[-0.8mm]{0.4mm}{3.4mm}\,\mu _{\omega }\,%
\rule[-0.8mm]{0.4mm}{3.4mm}\text{ }\Vert u\Vert _{\infty },  \label{aux-T1}
\end{equation}%
and therefore
\begin{equation*}
\Vert T\Vert _{\infty }\leq l_{0}\Vert f\Vert _{\infty }%
\rule[-0.8mm]{0.4mm}{3.4mm}\,\mu _{\omega }\,\rule[-0.8mm]{0.4mm}{3.4mm}.
\end{equation*}

For the nonlinear term we define $B:\mathcal{X}\rightarrow \mathcal{X}$ by
setting
\begin{equation*}
B(u):=H(b|u|^{p-1}u),~~~\forall \,u\in \mathcal{X}.  \label{aux-B}
\end{equation*}%
If $a_1,\,a_2 \in \mathbb{R}$ there holds
\begin{equation*}
\left\vert a_{1}|a_{1}|^{p-1}-a_{2}|a_{2}|^{p-1}\right\vert \leq
p|a_{1}-a_{2}|\left( |a_{1}|^{p-1}-|a_{2}|^{p-1}\right),
\end{equation*}%
and therefore it follows that
\begin{equation*}
\Vert b(\cdot)\left( u|u|^{p-1}-\tilde{u}|\tilde{u}|^{p-1}\right) \Vert
_{\infty }\leq \Vert b\Vert _{\infty }\Vert u-\tilde{u}\Vert _{\infty
}\left( \Vert u\Vert _{\infty }^{p-1}-\Vert \tilde{u}\Vert _{\infty
}^{p-1}\right).
\end{equation*}%
This inequality and the same argument used in (\ref{aux-T1}) imply that
\begin{equation}
\Vert B(u)-B(\tilde{u})\Vert _{\infty }\leq l_{0}p\Vert b\Vert _{\infty
}\Vert u-\tilde{u}\Vert _{\infty }\left( \Vert u\Vert _{\infty }^{p-1}-\Vert
\tilde{u}\Vert _{\infty }^{p-1}\right) ,  \label{est-nao-linear}
\end{equation}%
for any $u,\,\tilde{u}\in L^{\infty }(U)$.

All together, the above estimates enable us to solve the random equation %
\eqref{eqdif-est1} as follows.

\begin{proposition}
\label{ponto-fixo} Given $f,\,b,g\in L^{\infty }(U)$ and $\omega \in \Omega $%
, we consider the potential $V_{\omega }$ induced by the random measure $\mu
_{\omega }:=X(\omega )$. Let $l_{0}$ be the quantity introduced in %
\eqref{def-l0} and set
\begin{equation}
\tau _{\omega }:=l_{0}\Vert f\Vert _{\infty }\rule[-0.8mm]{0.4mm}{3.4mm}%
\,\mu _{\omega }\,\rule[-0.8mm]{0.4mm}{3.4mm}\text{ \ \ and}~~K:=l_{0}p\Vert
b\Vert _{\infty }.  \label{aux-tauK}
\end{equation}%
If $\varepsilon >0$ and $\omega \in \Omega $ are such that
\begin{equation}
0\leq \tau _{\omega }<1,~~~~~\frac{2^{p}K\varepsilon ^{p-1}}{(1-\tau
_{\omega })^{p-1}}+\tau _{\omega }<1,  \label{est-epsilon}
\end{equation}%
and $\Vert g\Vert _{\infty }\leq \varepsilon /l_{0}$, then the equation %
\eqref{eqdif-est1} has a unique integral solution
\begin{equation}
u_{\omega }=u(\cdot ,\omega )\in L^{\infty }(U)\text{ such that }\Vert
u_{\omega }\Vert _{\infty }\leq \frac{2\varepsilon }{1-\tau _{\omega }}.
\label{aux-sol1}
\end{equation}
\end{proposition}

\bigskip \noindent \textbf{Proof. }For each $\omega \in \Omega ,$ we
consider the closed ball
\begin{equation*}
\mathcal{B}_{\varepsilon }=\left\{ u\in L^{\infty }(U);\Vert u\Vert _{\infty
}\leq \frac{2\varepsilon }{(1-\tau _{\omega })}\right\}
\end{equation*}%
endowed with the metric $d(u,v):=\Vert u-v\Vert _{\infty }.$ We are going to
show that the map
\begin{equation}
\Phi (u):=H(g)+H(V_{\omega }u)+H(bu\left\vert u\right\vert
^{p-1})=H(g)+T(u)+B(u)  \label{aux-phi}
\end{equation}%
is a contraction on the complete metric space $(\mathcal{B}_{\varepsilon
},d).$ Using the estimates (\ref{est-h}), (\ref{aux-T1}), and (\ref%
{est-nao-linear}) with $\tilde{u}=0,$ we obtain
\begin{align}
\left\Vert \Phi (u)\right\Vert _{\infty }& \leq \Vert H(g)\Vert _{\infty
}+\Vert T(u)\Vert _{\infty }+\Vert B(u)\Vert _{\infty }  \notag \\[0.3cm]
& \leq l_{0}\left\Vert g\right\Vert _{\infty }+\tau _{\omega }\Vert u\Vert
_{\infty }+K\Vert u\Vert _{\infty }^{p}  \notag \\[0.3cm]
& \leq \varepsilon +\tau _{\omega }\frac{2\varepsilon }{1-\tau _{\omega }}+%
\frac{2^{p}K\varepsilon ^{p}}{(1-\tau _{\omega })^{p}}  \notag \\[0.3cm]
& =\left( 1+\tau _{\omega }+\frac{2^{p}K\varepsilon ^{p-1}}{(1-\tau _{\omega
})^{p-1}}\right) \frac{\varepsilon }{1-\tau _{\omega }}  \notag
\end{align}%
for all $u\in \mathcal{B}_{\varepsilon }$ and $\omega \in \Omega $. Hence,
it follows from (\ref{est-epsilon}) that
\begin{equation*}
\left\Vert \Phi (u)\right\Vert _{\infty }\leq \frac{2\varepsilon }{1-\tau
_{\omega }}.
\end{equation*}%
This shows that $\Phi $ maps $\mathcal{B}_{\varepsilon }$ into $\mathcal{B}%
_{\varepsilon }$.

For any $u,\widetilde{u}\in \mathcal{B}_{\varepsilon },$ it follows from (%
\ref{aux-T1}) and (\ref{est-nao-linear}) that
\begin{align}
\Vert \Phi (u)-\Phi (\widetilde{u})\Vert _{\infty }& =\Vert T(u-\widetilde{u}%
)\Vert _{\infty }+\Vert B(u)-B(\widetilde{u})\Vert _{\infty }  \notag \\%
[0.2cm]
& \leq \tau _{\omega }\Vert u-\widetilde{u}\Vert _{\infty }+K\Vert u-%
\widetilde{u}\Vert _{\infty }\left( \Vert u\Vert _{\infty }^{p-1}+\Vert
\widetilde{u}\Vert _{\infty }^{p-1}\right)  \notag \\[0.2cm]
& \leq \left( \tau _{\omega }+\frac{2^{p}K\varepsilon ^{p-1}}{(1-\tau
_{\omega })^{p-1}}\right) \Vert u-\widetilde{u}\Vert _{\infty }.  \notag
\label{aux-point2}
\end{align}%
Recalling (\ref{est-epsilon}), the above estimate implies that the map $\Phi
$ is a contraction. The Banach fixed point theorem assures that there is a
unique solution $u$ for the integral equation (\ref{eq-integral}) such that $%
\Vert u\Vert _{\infty }\leq (2\varepsilon)/(1-\tau _{\omega }).$ \fin

The next results are related to the randomness introduced by the random
potential $V$ and the existence and uniqueness of solutions for the problem %
\eqref{eqdif-est1}. Roughly speaking, we first obtain the probability of %
\eqref{eqdif-est1} having a solution given by the method discussed above. In
the sequel we study two important limit theorems in probability theory,
namely, the central limit theorem and the law of large numbers for a
sequence of random potentials.

\begin{theorem}
\label{teo-prob}Let $\nu $ be the probability measure induced on $\mathbb{R}$
by the random variable $\omega \mapsto \rule[-0.8mm]{0.4mm}{3.4mm}\,\mu
_{\omega }\rule[-0.8mm]{0.4mm}{3.4mm}\,\,.$ Let $g\in L^{\infty }(U)$ be
such that $\left\Vert g\right\Vert _{\infty }<\frac{1}{l_{0}}(\frac{1}{2^{p}K%
})^{\frac{1}{p-1}}$, where $K=l_{0}p\Vert b\Vert _{\infty }$. Choose $%
0<c_{0}<1$and set
\begin{equation*}
\varepsilon _{0}:=\left( \frac{(1-c_{0})^{p}}{2^{p}K}\right) ^{\frac{1}{p-1}%
}.
\end{equation*}%
Let $\mathcal{A}$ be the set of $\omega \in \Omega $ such that %
\eqref{eqdif-est1} has a unique solution $u(\cdot ,\omega )$ given by
Proposition \ref{ponto-fixo} with $\varepsilon =\varepsilon _{0}$. The set $%
\mathcal{A}$ is called the admissible one for the random variable $X.$

\begin{itemize}
\item[(i)] The set $\mathcal{A}$ is $\mathcal{F}$-measurable and the
probability of \eqref{eqdif-est1} having a solution is
\begin{equation*}
\mathbb{P(\mathcal{A})}=\nu \left( \left[ 0,\frac{1}{l_{0}\Vert f\Vert
_{\infty }}\right) \right) .
\end{equation*}

\item[(ii)] Let $u_{\omega },\,\tilde{u}_{\omega }$ be two solutions of %
\eqref{eqdif-est1} corresponding, respectively, to $\mu _{\omega },g,%
\mathcal{A}$ $\ $and $\tilde{\mu}_{\omega },\tilde{g},\widetilde{\mathcal{A}}
$. Assume that $\mathcal{A\cap }\widetilde{\mathcal{A}}\neq \varnothing $
and define, for $\omega \in \mathcal{A\cap }\widetilde{\mathcal{A}}$,
\begin{equation*}
\eta _{\omega }:=l_{0}\Vert f\Vert _{\infty }\max \{%
\rule[-0.8mm]{0.4mm}{3.4mm}\,\mu _{\omega }\,\rule[-0.8mm]{0.4mm}{3.4mm},%
\rule[-0.8mm]{0.4mm}{3.4mm}\,\widetilde{\mu }_{\omega }\,%
\rule[-0.8mm]{0.4mm}{3.4mm}\}.
\end{equation*}%
We have that
\begin{equation}
\Vert u(\cdot ,\omega )-\tilde{u}(\cdot ,\omega )\Vert _{\infty }\leq \frac{%
l_{0}\left( \Vert g-\tilde{g}\Vert _{\infty }+\dfrac{2\varepsilon _{0}}{%
1-\eta _{\omega }}\Vert f\Vert _{\infty }\rule[-0.8mm]{0.4mm}{3.4mm}\,\mu
_{\omega }-\tilde{\mu}_{\omega }\,\rule[-0.8mm]{0.4mm}{3.4mm}\right) }{%
1-\eta _{\omega }-\dfrac{2^{p}K\varepsilon _{0}^{p-1}}{(1-\eta _{\omega
})^{p-1}}}  \label{lip}
\end{equation}%
for all $\omega \in \mathcal{A\cap }\widetilde{\mathcal{A}}$.

\item[(iii)] The map $\ \mathcal{U}:\mathcal{A}\rightarrow L^{\infty }(U)$
given by $\mathcal{U}(\omega ):=u(\cdot ,\omega )$ is a random variable and
there holds
\begin{equation}
\Vert u(\cdot ,\omega )\Vert _{\infty }\leq \frac{2\varepsilon _{0}}{1-\tau
_{\omega }}=2\varepsilon _{0}\sum_{j=0}^{\infty }\tau _{\omega }^{j},
\label{aux-serie}
\end{equation}%
for all $\omega \in \mathcal{A}.$
\end{itemize}
\end{theorem}

\bigskip

\noindent \textbf{Proof.} We first notice that the choice of $\varepsilon
_{0}$ implies that $\Vert g\Vert _{\infty }\leq \varepsilon _{0}/l_{0}.$
Moreover, $\omega \in \mathcal{A}$ if only if $\tau _{\omega }=l_{0}\Vert
f\Vert _{\infty }\rule[-0.8mm]{0.4mm}{3.4mm}\,\mu _{\omega }\,%
\rule[-0.8mm]{0.4mm}{3.4mm}$ verifies (\ref{est-epsilon}) with $\varepsilon
=\varepsilon _{0}.$ Then, if $Y(\omega )=\rule[-0.8mm]{0.4mm}{3.4mm}%
\,X(\omega )\rule[-0.8mm]{0.4mm}{3.4mm}\,=\rule[-0.8mm]{0.4mm}{3.4mm}\,\,\mu
_{\omega }\,\rule[-0.8mm]{0.4mm}{3.4mm},$ it follows that $\mathcal{A}%
=\left\{ Y\in \left[ 0,\frac{1}{l_{0}\Vert f\Vert _{\infty }}\right)
\right\} $ is measurable and

\begin{equation*}
\begin{array}{lcl}
\mathbb{P}(\mathcal{A}) & = & \mathbb{P}\left( Y\in \left[ 0,\dfrac{1}{%
l_{0}\Vert f\Vert _{\infty }}\right) \right) =\mathbb{P}_{Y}\left( \left[ 0,%
\dfrac{1}{l_{0}\Vert f\Vert _{\infty }}\right) \right) \vspace{0.2cm} \\
& = & \nu \left( \left[ 0,\dfrac{1}{l_{0}\Vert f\Vert _{\infty }}\right)
\right) . \label{aux-prob1}%
\end{array}%
\end{equation*}%
This establishes (i).

Now we deal with item (ii). Firstly, observe that $\eta _{\omega }=\max
\{\tau _{\omega },\widetilde{\tau }_{\omega }\},$ where
\begin{equation*}
\tau _{\omega }=l_{0}\Vert f\Vert _{\infty }\rule[-0.8mm]{0.4mm}{3.4mm}\mu
_{\omega }\,\rule[-0.8mm]{0.4mm}{3.4mm}\text{ and }\widetilde{\tau }_{\omega
}=l_{0}\Vert f\Vert _{\infty }\rule[-0.8mm]{0.4mm}{3.4mm}\tilde{\mu}_{\omega
}\,\rule[-0.8mm]{0.4mm}{3.4mm}.
\end{equation*}%
Subtracting the integral equations verified by $u_{\omega }$ and $\,\tilde{u}%
_{\omega },$ and afterwards computing $\Vert \cdot \Vert _{\infty }$, we
obtain
\begin{equation*}
\begin{array}{lcl}
\left\Vert u_{\omega }-\tilde{u}_{\omega }\right\Vert _{\infty } & \leq &
\left\Vert H(g-\tilde{g})\right\Vert _{\infty }+\left\Vert H(V_{\omega }(u-%
\tilde{u}_{\omega }))\right\Vert _{\infty }\vspace{0.2cm} \\
&  & +\Vert H((V_{\omega }-\widetilde{V}_{\omega })\tilde{u}_{\omega })\Vert
_{\infty }\vspace{0.2cm} \\
&  & +\left\Vert H(b\left( u_{\omega }|u_{\omega }|^{p-1}-\tilde{u}_{\omega
}|\tilde{u}_{\omega }|^{p-1}\right) )\right\Vert _{\infty }\vspace{0.2cm} \\
& \leq & l_{0}\Vert g-\tilde{g}\Vert _{\infty }+l_{0}\Vert f\Vert _{\infty }%
\rule[-0.8mm]{0.4mm}{3.4mm}\,\mu _{\omega }\,\rule[-0.8mm]{0.4mm}{3.4mm}%
\Vert u_{\omega }-\tilde{u}_{\omega }\Vert _{\infty }\vspace{0.2cm} \\
&  & +l_{0}\Vert f\Vert _{\infty }\rule[-0.8mm]{0.4mm}{3.4mm}\,\mu _{\omega
}-\tilde{\mu}_{\omega }\,\rule[-0.8mm]{0.4mm}{3.4mm}\Vert \tilde{u}_{\omega
}\Vert _{\infty }\vspace{0.2cm} \\
&  & +l_{0}p\Vert b\Vert _{\infty }\Vert u_{\omega }-\tilde{u}_{\omega
}\Vert _{\infty }(\Vert u_{\omega }\Vert _{\infty }^{p-1}-\Vert \tilde{u}%
_{\omega }\Vert _{\infty }^{p-1}).%
\end{array}%
\end{equation*}%
It follows from (\ref{aux-sol1}) that
\begin{equation*}
\Vert u_{\omega }\Vert _{\infty }\leq \frac{2\varepsilon _{0}}{1-\tau
_{\omega }}\leq \frac{2\varepsilon _{0}}{1-\eta _{\omega }}\text{ and }\Vert
\tilde{u}\Vert _{\infty }\leq \frac{2\varepsilon _{0}}{1-\widetilde{\tau }%
_{\omega }}\leq \frac{2\varepsilon _{0}}{1-\eta _{\omega }}.
\end{equation*}%
The two above expressions give us
\begin{eqnarray*}
\left\Vert u_{\omega }-\tilde{u}_{\omega }\right\Vert _{\infty } &\leq
&l_{0}\Vert g-\tilde{g}\Vert _{\infty }+l_{0}\Vert f\Vert _{\infty }%
\rule[-0.8mm]{0.4mm}{3.4mm}\,\mu _{\omega }\,\rule[-0.8mm]{0.4mm}{3.4mm}%
\left\Vert u_{\omega }-\tilde{u}_{\omega }\right\Vert _{\infty } \\[0.3cm]
&&+\,l_{0}\frac{2\varepsilon _{0}}{1-\eta _{\omega }}\Vert f\Vert _{\infty }%
\rule[-0.8mm]{0.4mm}{3.4mm}\,\mu _{\omega }-\tilde{\mu}_{\omega }\,%
\rule[-0.8mm]{0.4mm}{3.4mm}+\frac{2^{p}K\varepsilon _{0}^{p-1}}{(1-\eta
_{\omega })^{p-1}}\Vert u_{\omega }-\tilde{u}_{\omega }\Vert _{\infty } \\%
[0.4cm]
&=&l_{0}\Vert g-\tilde{g}\Vert _{\infty }+l_{0}\frac{2\varepsilon _{0}}{%
1-\eta _{\omega }}\Vert f\Vert _{\infty }\rule[-0.8mm]{0.4mm}{3.4mm}\,\mu
_{\omega }-\tilde{\mu}_{\omega }\,\rule[-0.8mm]{0.4mm}{3.4mm} \\[0.3cm]
&&+\,\left[ \eta _{\omega }+\frac{2^{p}K\varepsilon _{0}^{p-1}}{(1-\eta
_{\omega })^{p-1}}\right] \left\Vert u_{\omega }-\tilde{u}_{\omega
}\right\Vert _{\infty },
\end{eqnarray*}%
which yields (\ref{lip}).

Taking $\mu _{\omega },\tilde{\mu}_{\omega }$ independent of $\omega ,$ i.e.
$\mu _{\omega }=\mu $ and $\tilde{\mu}_{\omega }=\tilde{\mu},$ for all $%
\omega \in \Omega ,$ we see from (\ref{aux-tauK}) and (\ref{lip}) that the
data-map solution $\mathcal{L}(\mu ,g)=u$ is continuous from
\begin{equation}
\hspace*{-0.4cm}\left\{ (\mu ,g)\in \mathcal{M}(U)\times L^{\infty }(U);%
\text{ }\rule[-0.8mm]{0.4mm}{3.4mm}\,\mu \,\rule[-0.8mm]{0.4mm}{3.4mm}<\frac{%
1}{l_{0}\Vert f\Vert _{\infty }},\left\Vert g\right\Vert _{\infty }<\frac{1}{%
l_{0}}\left( \frac{1}{2^{p}K}\right) ^{\frac{1}{p-1}}\right\} \text{to }%
L^{\infty }(U),  \label{aux-proof2}
\end{equation}%
where $u$ is the deterministic solution of (\ref{eqdif-est1}) corresponding
to the data $(\mu ,g).$ From this, and because $X|_{\mathcal{A}}$ given by $%
X(\omega )=\mu _{\omega }$ is measurable, it follows that the composition $%
\mathcal{U}(\omega )=\mathcal{L}(\mu _{\omega },g)=\mathcal{L}(X(\omega
),g)\,$\ from $\mathcal{A}$ to $L^{\infty }(U)$ is measurable.

In view of the series $\frac{1}{1-z}=\sum_{j=0}^{\infty }z^{j}$ for $%
\left\vert z\right\vert <1,$ we finish by observing that (\ref{aux-serie})
follows at once from (\ref{aux-sol1}) with $\varepsilon =\varepsilon _{0}$
and $\omega \in \mathcal{A}.$\fin

\begin{remark}
\label{Rem1} Here we give examples of random potentials for which there
exists solution almost surely in $\Omega $. The first setting occurs if we
suppose that the measure $\nu $ $\,$has compact support contained in the
interval $[0,a]$, with $a<\frac{1}{l_{0}\Vert f\Vert _{\infty }}$. In this
case it follows from the first item of the above theorem that $\mathbb{P(%
\mathcal{A})}=1$, i.e., the solution exists almost surely in $\Omega .$
Secondly, we take $\{\mu _{j}\}_{j\in \mathbb{N}}$ a sequence in $\mathcal{M}%
(U)$ and let $\{a_{j}(\omega )\}_{j\in \mathbb{N}}$ be a sequence of random
variables from $\Omega $ to $\mathbb{R}.$ Consider the random variable $\mu
_{\omega }$ defined by
\begin{equation*}
\mu _{\omega }=\sum_{j=1}^{\infty }a_{j}(\omega )\mu _{j}.
\end{equation*}%
For $q>1,$ suppose that
\begin{equation*}
|a_{j}(\omega )|<\frac{(\sum_{k=1}^{\infty }\frac{1}{k^{q}})^{-1}}{l_{0}%
\rule[-0.8mm]{0.4mm}{3.4mm}\,\mu _{j}\rule[-0.8mm]{0.4mm}{3.4mm}\,\Vert
f\Vert _{\infty }}\cdot \frac{1}{j^{q}}\text{ a.s. in }\Omega ,
\end{equation*}%
for all $j\in \mathbb{N}.$ Then%
\begin{equation*}
\rule[-0.8mm]{0.4mm}{3.4mm}\,\mu _{\omega }\rule[-0.8mm]{0.4mm}{3.4mm}\,\leq
\sum_{j=1}^{\infty }|a_{j}(\omega )|\rule[-0.8mm]{0.4mm}{3.4mm}\,\mu _{j}%
\rule[-0.8mm]{0.4mm}{3.4mm}\,<\frac{1}{l_{0}\Vert f\Vert _{\infty }}\text{
a.s. in }\Omega ,
\end{equation*}%
and Theorem \ref{teo-prob} assures that there is an integral solution for (%
\ref{eqdif-est1}) a.s. in $\Omega .$
\end{remark}

In the sequel we show how the Borel-Cantelli's Lemma can be used to give a
sufficient condition for the existence of solution a.s. in $\Omega $.

\begin{theorem}
\label{borel-canteli} Let $\{\mu _{j}\}_{j\in \mathbb{N}}$ be a sequence in $%
\mathcal{M}(U)$ and let $\{a_{j}(\omega )\}_{j\in \mathbb{N}}$ be a sequence
of random variables from $\Omega $ to $\mathbb{R}.$ Assume that the
following series is convergent in $\mathcal{M}(U)$
\begin{equation*}
\mu _{\omega }=\sum_{j=1}^{\infty }a_{j}(\omega )\mu _{j}.
\end{equation*}%
For any $k\in \mathbb{N}$ define
\begin{equation*}
S_{k}(\omega )=\sum_{j=1}^{k}a_{j}(\omega )\mu _{j}
\end{equation*}%
and $L_{k}=\{\omega \in \Omega :\rule[-0.8mm]{0.4mm}{3.4mm}\,S_{k}%
\rule[-0.8mm]{0.4mm}{3.4mm}\,\geq \tilde{c}\}$, with $0<\tilde{c}%
<1/(l_{0}\Vert f\Vert _{\infty })$. If
\begin{equation*}
\sum_{k=1}^{\infty }\mathbb{P}(L_{k})<\infty
\end{equation*}%
then there is an integral solution for \eqref{eqdif-est1} almost surely in $%
\Omega $.
\end{theorem}

\noindent \textbf{Proof. \ }By the Borel-Cantelli's Lemma we get that $%
\mathbb{P}(\limsup L_{k})=0$, that is,
\begin{equation*}
\mathbb{P}\left( \cup _{j=1}^{\infty }\cap _{k=j}^{\infty }\{%
\rule[-0.8mm]{0.4mm}{3.4mm}\,S_{k}\,\rule[-0.8mm]{0.4mm}{3.4mm}<\tilde{c}%
\}\right) =1
\end{equation*}%
It follows that, for almost sure $\omega ,$ there is $j_{0}=$ $j_{0}(\omega )
$ such that for all $j>j_{0}$, we have
\begin{equation*}
\rule[-0.8mm]{0.4mm}{3.4mm}\,S_{k}(\omega )\,\rule[-0.8mm]{0.4mm}{3.4mm}<%
\tilde{c}.
\end{equation*}%
Therefore by taking the limit when $k$ goes to infinity, we obtain
\begin{equation*}
\rule[-0.8mm]{0.4mm}{3.4mm}\,\mu _{\omega }\,\rule[-0.8mm]{0.4mm}{3.4mm}%
=\lim_{j\rightarrow \infty }\rule[-0.8mm]{0.4mm}{3.4mm}\,S_{j}(\omega )\,%
\rule[-0.8mm]{0.4mm}{3.4mm}\leq \tilde{c}<\frac{1}{l_{0}\Vert f\Vert
_{\infty }}\quad \mbox{a.s. in }\Omega .
\end{equation*}%
This inequality and Theorem \ref{teo-prob} imply that there is an integral
solution $u(x,\omega )$ for \eqref{eqdif-est1} almost surely in $\Omega $. %
\fin

\bigskip

A straightforward calculation shows that in general $\mathbb{E}_{\Omega
}(u(x,\omega ))$ does not satisfies the equation \eqref{eqdif-est1}, even if
we replace the random potential by its mean. However, we are able to obtain
some information on the average and moments of the random solution $%
u_{\omega }$ previously obtained. It is worthwhile to mention
that, when dealing with the random variable $\omega \mapsto
u_{\omega }$, the expectation has to be understood in the Bochner
sense (see Section 2). Note also that a solution $u_{\omega }\in
L^{\infty }(U)$ for (\ref{eq-integral}) in fact belongs to the
separable subspace $C(\overline{U}).$

\begin{theorem}
\label{teo-moment}Under hypotheses of Theorem \ref{teo-prob} let us denote
by $u_{\omega }(x)=u(x,\omega )\in \mathcal{A}$ the solution of %
\eqref{eqdif-est1}. Let $m\in \mathbb{N}$ and suppose that
\begin{equation}
\sum_{j=1}^{\infty }\frac{(m+j-1)!}{(m-1)!j!}(l_{0}\Vert f\Vert _{\infty
})^{j}\ \mathbb{E}_{\mathcal{A}}[\ \rule[-0.8mm]{0.4mm}{3.4mm}\,\mu _{\omega
}\,\rule[-0.8mm]{0.4mm}{3.4mm}^{\,j}\,]<+\infty .  \label{hip-moment2}
\end{equation}%
Then $\mathbb{E}_{\mathcal{A}}[|u|^{m}(x,\omega )]\in L^{\infty }(U)$ and
\begin{equation}
\mathbb{E}_{\mathcal{A}}\left[ \left\Vert |u|^{m}(\cdot ,\omega )\right\Vert
_{L^{\infty }(U)}\right] <\infty .  \label{moment}
\end{equation}
In particular, $\mathbb{E}_{\mathcal{A}}[u(x,\omega )]\in L^{\infty }(U).$
\end{theorem}

\noindent\textbf{Proof.} It follows from (\ref{aux-serie}) that

\begin{equation}
\Vert |u|^{m}(\cdot ,\omega )\Vert _{L^{\infty }(U)}\leq \Vert u(\cdot
,\omega )\Vert _{L^{\infty }(U)}^{m}\leq \frac{(2\varepsilon _{0})^{m}}{%
(1-\tau _{\omega })^{m}}.  \label{aux-moment1}
\end{equation}%
Recalling that $\tau _{\omega }=l_{0}\Vert f\Vert _{\infty }%
\rule[-0.8mm]{0.4mm}{3.4mm}\,\mu _{\omega }\,\rule[-0.8mm]{0.4mm}{3.4mm}$
and computing $\mathbb{E}_{\mathcal{A}}$ in (\ref{aux-moment1}), we obtain
\begin{eqnarray}
\left\Vert \mathbb{E}_{\mathcal{A}}\left[ |u|^{m}(x,\omega )\right]
\right\Vert _{L^{\infty }(U)} &\leq &\mathbb{E}_{\mathcal{A}}\left[
\left\Vert |u|^{m}(x,\omega )\right\Vert _{L^{\infty }(U)}\right]  \notag \\
&\leq &(2\varepsilon _{0})^{m}\mathbb{E}_{\mathcal{A}}\left[ \left(
1+\sum_{j=1}^{\infty }\frac{(m+j-1)!}{(m-1)!j!}\tau _{\omega }^{j}\right) %
\right]  \notag
\end{eqnarray}%
By using the linearity of the expectation and definition of $\tau _{\omega }$
we get the following upper bound for the right hand side above
\begin{equation*}
(2\varepsilon _{0})^{m}+(2\varepsilon _{0})^{m}\sum_{j=1}^{\infty }\frac{%
(m+j-1)!}{(m-1)!j!}\left( l_{0}\Vert f\Vert _{\infty }\right) ^{j}\mathbb{E}%
_{\mathcal{A}}\left[ \ \rule[-0.8mm]{0.4mm}{3.4mm}\,\mu _{\omega }\,%
\rule[-0.8mm]{0.4mm}{3.4mm}^{\,j}\ \right] ,
\end{equation*}%
which is finite due to (\ref{hip-moment2}). The last assertion of the
statement follows from (\ref{moment}) with $m=1$ and the easy estimate
\begin{equation*}
\left\Vert \mathbb{E}_{\mathcal{A}}\left[ u(x,\omega )\right] \right\Vert
_{L^{\infty }(U)}\leq \mathbb{E}_{\mathcal{A}}\left[ \left\Vert |u|(x,\omega
)\right\Vert _{L^{\infty }(U)}\right] .
\end{equation*}%
\fin

\subsection{Classical Probability Limit Theorems}

We start this section by recalling basic background concerning to some main
limit theorems in probability. A real-valued random variable $X:\Omega\to%
\mathbb{R}$ in a probability space $(\Omega,\mathcal{F},\mathbb{P})$ has
standard normal distribution, notation $X\sim N(0,1)$, if for all $x \in
\mathbb{R}$ its cumulative distribution function verifies
\begin{equation*}
\mathbb{P}(X\leq x) = \frac{1}{\sqrt{2\pi}} \int_{-\infty}^x e^{-\frac{1}{2}%
t^2} dt.
\end{equation*}
A sequence of real-valued random variable $\{Y_j\}_{j\in\mathbb{N}}$ in a
probability space $(\Omega,\mathcal{F},\mathbb{P})$ is said to converge in
distribution to a standard normal random variable, notation $Y_j\to N(0,1)$,
if for all $x\in\mathbb{R}$ we have
\begin{equation*}
\lim_{j\to\infty} \mathbb{P}(Y_j\leq x)=\frac{1}{\sqrt{2\pi}}
\int_{-\infty}^x e^{-\frac{1}{2}t^2} dt.
\end{equation*}

In the sequel we show versions of the central limit theorem and a weak law
of large numbers for the random $L^{\infty}(U)$-solutions obtained in
Section 2.

\begin{theorem}
Let $\{X_{j}\}_{j\in \mathbb{N}}$ be an independent identically distributed
(i.i.d.) sequence of random variables $X_{j}:\Omega \rightarrow \mathcal{M}%
(U)$. Assume that the admissible set $\mathcal{A}_{j}=\Omega $ for all $j$,
and let $u_{j}(\cdot ,\omega )\in L^{\infty }(U)$ be the solution given by
Theorem \ref{teo-prob} with respect to $X_{j}(\omega )=\mu _{\omega ,j}$ and
$g.$ We have that $\{Z_{j}\}_{j\in \mathbb{N}}$ given by $Z_{j}(\omega
):=\Vert u_{j}(\cdot ,\omega )\Vert _{\infty }$ is a i.i.d. sequence of
random variables, and if $m=\mathbb{E}[\Vert u_{j}(\cdot ,\omega )\Vert
_{\infty }]<\infty $ and $\sigma ^{2}:=\text{Var}\,Z_{j}<\infty $ then
following holds as $k\rightarrow +\infty $
\begin{equation*}
\sum_{j=1}^{k}\frac{(Z_{j}-m)}{\sigma \sqrt{k}}\rightarrow N(0,1).
\end{equation*}
\end{theorem}

\noindent \textbf{Proof.} Recall the data-solution map $\mathcal{L}(\mu ,g)$
defined in the proof of Theorem \ref{teo-prob} (see (\ref{aux-proof2})).
Fixed $g$ such that $\left\Vert g\right\Vert _{\infty }<\frac{1}{l_{0}}(%
\frac{1}{2^{p}K})^{\frac{1}{p-1}},$ consider
\begin{equation}
S_{g}(\mu )=\mathcal{L}(\mu ,g)  \label{map-Sg}
\end{equation}%
defined from $D$ to $L^{\infty }(U),$ where $D=\left\{ \mu \in \mathcal{M}%
(U):\rule[-0.8mm]{0.4mm}{3.4mm}\,\mu \,\rule[-0.8mm]{0.4mm}{3.4mm}<\frac{1}{%
l_{0}\Vert f\Vert _{\infty }}\right\} .$ Since $\Vert \cdot \Vert _{\infty }$
is continuous from $L^{\infty }(U)$ to $\mathbb{R}$ and
\begin{equation*}
Z_{j}(\omega )=\Vert u_{j}(\cdot ,\omega )\Vert _{\infty }=\Vert S_{g}\circ
X_{j}(\omega )\Vert _{\infty },
\end{equation*}%
we get that $\{Z_{j}\}_{j\in \mathbb{N}}$ is a i.i.d. sequence. The
convergence stated in the theorem follows from the central limit theorem. %
\fin

\begin{theorem}
Let $\{X_{j}\}_{j\in \mathbb{N}}$ be an independent sequence of random
variables $X_{j}:\Omega \rightarrow \mathcal{M}(U)$. Assume that the
admissible set $\mathcal{A}_{j}=\Omega $ for all $j$, and let $u_{j}(\cdot
,\omega )\in L^{\infty }(U)$ be the solution given by Theorem \ref{teo-prob}
with respect to $X_{j}(\omega )=\mu _{\omega ,j}$ and $g.$ If $%
X_{j}\rightarrow X$ a.s. and
\begin{equation}
L=\sup_{j\in \mathbb{N}}\left( \mathrm{ess}\sup_{\omega \in \Omega }%
\rule[-0.8mm]{0.4mm}{3.4mm}\,\mu _{\omega ,j}\,\rule[-0.8mm]{0.4mm}{3.4mm}%
\right) <\frac{1}{l_{0}\Vert f\Vert _{\infty }},  \label{hip-law}
\end{equation}%
then%
\begin{equation}
\sum_{j=1}^{k}\frac{u_{j}(x,\omega )-\mathbb{E}_{\Omega }[u_{j}(x,\omega )]}{%
k}\rightarrow 0  \label{conv1}
\end{equation}%
and
\begin{equation}
\sum_{j=1}^{k}\frac{\Vert u_{j}(\cdot ,\omega )\Vert _{\infty }-\mathbb{E}%
_{\Omega }[\Vert u_{j}(\cdot ,\omega )\Vert _{\infty }]}{k}\rightarrow 0,
\label{conv2}
\end{equation}%
when $k\rightarrow \infty $, where the convergence in (\ref{conv1}) and (\ref%
{conv2}) are in probability sense (see (\ref{def-conv-prob})).
\end{theorem}

\noindent \textbf{Proof.} Notice that $X_{j}\rightarrow X$ a.s. is
equivalent to $\mu _{\omega ,j}\rightarrow \mu _{\omega }=X(\omega )$ in $%
\mathcal{M}(U)$ almost surely. From this and the continuity of data-solution
map $\mathcal{L}(\cdot ,\cdot )$ (see (\ref{aux-proof2})), it follows that
\begin{equation*}
\Vert u_{j}(\cdot ,\omega )-u(\cdot ,\omega )\Vert _{\infty }=\Vert \mathcal{%
L}(\mu _{\omega ,j},g)-\mathcal{L}(\mu ,g)\Vert _{\infty }\rightarrow 0,
\end{equation*}%
when $j\rightarrow \infty $. Recalling (\ref{aux-serie}) and afterwards
using (\ref{hip-law}), we obtain
\begin{eqnarray}
\Vert u_{j}(\cdot ,\omega )\Vert _{\infty } &\leq &\frac{2\varepsilon _{0}}{%
1-l_{0}\Vert f\Vert _{\infty }(\mathrm{ess}\sup_{\omega \in \Omega }%
\rule[-0.8mm]{0.4mm}{3.4mm}\,\mu _{\omega },_{j}\,\rule[-0.8mm]{0.4mm}{3.4mm}%
)}  \notag \\
&\leq &\frac{2\varepsilon _{0}}{1-L}=Q_{0},\text{ a.s. in }\Omega .
\label{aux-proof-law} \\
&&  \notag
\end{eqnarray}%
Since $X_{j}$'s are independent, it follows that $\{Y_{j}\}_{j\in \mathbb{N}%
} $ defined by $Y_{j}=\left\Vert u_{j}(\cdot ,\omega )\right\Vert _{\infty
}=\Vert S_{g}\circ X_{j}(\omega )\Vert _{\infty }$ are also independent,
where $S_{g}$ is as in (\ref{map-Sg}). So, from Chebyshev's inequality and
the independence of $\{Y_{j}\}_{j\in \mathbb{N}}$, we have that
\begin{eqnarray*}
&&\hspace*{-2.5cm}\mathbb{P}\left( \left\vert k^{-1}\sum_{j=1}^{k}(\Vert
u_{j}(\cdot ,\omega )\Vert _{\infty }-\mathbb{E}_{\Omega }[\Vert u_{j}(\cdot
,\omega )\Vert _{\infty }])\right\vert \geq \delta \right) \\
&\leq &\frac{1}{(k\delta )^{2}}\mathbb{E}_{\Omega }\left[ \left\vert
\sum_{j=1}^{k}\left( \Vert u_{j}(\cdot ,\omega )\Vert _{\infty }-\mathbb{E}%
_{\Omega }[\Vert u_{j}(\cdot ,\omega )\Vert _{\infty }]\text{ }\right)
\right\vert ^{2}\right] \\
&=&\frac{1}{(k\delta )^{2}}\sum_{j=1}^{k}\mathbb{E}_{\Omega }\left[
\left\vert \left( \Vert u_{j}(\cdot ,\omega )\Vert _{\infty }-\mathbb{E}%
_{\Omega }[\Vert u_{j}(\cdot ,\omega )\Vert _{\infty }]\text{ }\right)
\right\vert ^{2}\right] \\
&\leq &\frac{1}{(k\delta )^{2}}\sum_{j=1}^{k}\mathbb{E}_{\Omega }\left[
\left\vert 2Q_{0}\right\vert ^{2}\right] \leq \frac{4Q_{0}^{2}}{\delta ^{2}}%
\frac{1}{k},
\end{eqnarray*}%
where we have used (\ref{aux-proof-law}). Letting $k\rightarrow +\infty $ in
the above expression we get (\ref{conv2}). The convergence (\ref{conv1}) can
be proved with similar arguments. \fin

\section*{Acknowledgments}

L.C.F. Ferreira was supported by FAPESP-SP and CNPq, Brazil. M. Furtado was
supported by CNPq, Brazil.

\end{document}